\documentclass[12pt,reqno]{article}

\usepackage[usenames]{color}
\usepackage{amssymb}
\usepackage{amsmath}
\usepackage{amsthm}
\usepackage{amsfonts}
\usepackage{amscd}
\usepackage{graphicx}
\usepackage[colorlinks=true,
linkcolor=webgreen,
filecolor=webbrown,
citecolor=webgreen]{hyperref}

\definecolor{webgreen}{rgb}{0,.5,0}
\definecolor{webbrown}{rgb}{.6,0,0}

\usepackage{color}
\usepackage{fullpage}
\usepackage{float}
\usepackage{graphics}
\usepackage{latexsym}

\begin{document}

\theoremstyle{plain}
\newtheorem{theorem}{Theorem}
\newtheorem{corollary}[theorem]{Corollary}
\newtheorem{proposition}{Proposition}
\newtheorem{lemma}{Lemma}
\newtheorem{example}{Examples}
\newtheorem{remark}{Remark}

\begin{center}
%{\large\bf
%	Families of Ramanujan-like Series for $1/\pi$ 
%Derived from Newton's Binomial Theorem
%}
%\vskip .5cm

\large\bf
A Master Identity for Forsyth-Glaisher-Type Series for $1/\pi$

\vskip 1.cm

{\large
Kunle Adegoke \\
Department of Physics and Engineering Physics \\ Obafemi Awolowo University, 220005 Ile-Ife, Nigeria \\
\href{mailto:adegoke00@gmail.com}{\tt adegoke00@gmail.com}
}

\end{center}

\vskip .15 in

\begin{abstract}
We derive a novel three-parameter master identity for Forsyth-Glaisher-type series for $1/\pi$. By unifying disparate results from the literature into a single overarching framework, we demonstrate that Glaisher's original series and numerous known variations emerge naturally as special cases through systematic choices of the parameters. This master identity simplifies the derivation of existing formulas, provides a mechanism for generating an infinite number of new series for $1/\pi$, and is further utilized to derive associated series involving harmonic and odd harmonic numbers.
\end{abstract}

\noindent 2010 {\it Mathematics Subject Classification}: Primary 33C05; Secondary 11Y60.

\noindent \emph{Keywords:} binomial coefficient; harmonic number, odd harmonic number, hypergeometric series, Forsyth-Glaisher-type series, Ramanujan-type series.

\bigskip

\section{Introduction}
Forsyth-Glaisher-type series for $1/\pi$ are a class of infinite series that typically feature central binomial coefficients $\binom {2k}k$ and rational factors. Historically, they were derived using Legendre polynomials, Fourier-Legendre expansions or elliptic functions.

A generic Forsyth-Glaisher-type series has the following form:
\begin{equation*}
\sum_{k = 0}^\infty  {\frac{{\binom{{2k}}{k}^2 }}{{P(k)c^k }}}  = \frac{A}{\pi },
\end{equation*}
where $P(k)$ is a polynomial in $k$ or $2k$ with real or fractional roots such as $(2k-1)^2$, $(k+1)^2$, or $(2k-3)(2k-1)$; $c$ is a geometric shifting base usually restricted to classical powers of $2$; and $A$ is a constant. 

Forsyth-Glaisher-type series represent a vital missing link in the history of $\pi$ formulas, bridging the gap between classical 19th-century orthogonal polynomials and the hyper-fast converging Ramanujan-like
hypergeometric series. Historically overshadowing simple Taylor expansions, these series rely on squaring
operations of binomial blocks or exploiting elliptic integral transforms to produce rational approximations for $1/\pi$.

Quintessential examples of Forsyth-Glaisher-type series are the original Forsyth's series derived by A.~R.~Forsyth in 1883~\cite{forsyth} using Legendre polynomials, namely,
\begin{equation}\label{forsyth}
	\sum_{k=0}^{\infty}
	\frac{\binom{2k}{k}^{2}}{2^{4k}(2k-1)^2}
	=\frac{4}{\pi},
\end{equation}
and the original Glaisher's series derived by J.~W.~L.~Glaisher in 1905~\cite{glaisher}, using the expansion of complete elliptic integrals:
\begin{equation}\label{glaisher}
	\sum_{k=0}^{\infty}
	\frac{\binom{2k}{k}^{2}}{2^{4k}(k+1)}
	=\frac{4}{\pi}.
\end{equation}

Let $\left( x \right)_k  = \prod_{j = 0}^{k - 1} {\left( {x + j} \right)} $, $x$ a complex number and $k$ a non-negative integer. Since
\begin{equation*}
\binom{2k}k=2^{2k}\,\frac{\left(\frac12\right)_k}{k!},
\end{equation*}
in terms of~$(x)_j$, the above prototype Forsyth-Glaisher-type series has the form:
\begin{equation*}
\sum_{k = 0}^\infty  {\frac{{\left( {\frac{1}{2}} \right)_k^2 }}{{k!^2}}\,\frac1{P(k)}}  = \frac{A}{\pi },
\end{equation*}
where $A$ and $P(k)$ are as previously defined.

Both~\eqref{forsyth} and~\eqref{glaisher} are special cases of the following master identity
\begin{align}\label{master1}
&\sum_{k = 0}^\infty  {\frac{{\left( {\frac{1}{2}} \right)_k^2 }}{{k!^2 }}\,\frac{1}{{\prod_{j = 1}^{\,\ell} {\left( {2k - 2j + 1} \right)} \prod_{j = 1}^m {\left( {2k - 2j + 1} \right)} \prod_{j = 1}^n {\left(k + j\right)} }}}\nonumber\\
&\qquad  = \frac{{\left( {\,\ell + m + n - 1} \right)!}}{{\left( {1/2} \right)_{\,\ell}^2 \left( {1/2} \right)_m^2 \left( {1/2} \right)_n^2 }}\,\frac{{( - 1)^{\,\ell + m} }}{{2^{\,\ell + m} }}\,\frac{{\binom{{\,\ell + n}}{n}\binom{{m + n}}{n}}}{{\binom{{2\,\ell + 2n}}{{2n}}\binom{{2m + 2n}}{{2n}}}}\,\frac{1}{\pi },
\end{align}
to be derived in this paper. Identity~\eqref{master1} is valid for non-negative integers $\ell$, $m$, and $n$ such that $\ell+m+n$ is a positive integer. The master identity subsumes many Glaisher series for $1/\pi$ enumerated in the paper by Nimbran~\cite{nimbran15} as well as the new series for $1/\pi$ derived therein.

Let $H_j$ and $O_j$ be the $j^{th}$ standard and odd harmonic number, respectively.

We will also derive the harmonic number analogue and the odd harmonic number analogue of~\eqref{master1}, namely,
\begin{align}\label{master2}
&\sum_{k = 0}^\infty  {\frac{{\left( {\frac{1}{2}} \right)_k^2 }}{{k!^2 }}\,\frac{H_{k+n}}{{\prod_{j = 1}^{\,\ell} {\left( {2k - 2j + 1} \right)} \prod_{j = 1}^m {\left( {2k - 2j + 1} \right)} \prod_{j = 1}^n {\left(k + j\right)} }}}\nonumber\\
&\qquad  = \frac{{\left( {\,\ell + m + n - 1} \right)!}}{{\left( {1/2} \right)_{\,\ell}^2 \left( {1/2} \right)_m^2 \left( {1/2} \right)_n^2 }}\,\frac{{( - 1)^{\,\ell + m} }}{{2^{\,\ell + m} }}\,\frac{{\binom{{\,\ell + n}}{n}\binom{{m + n}}{n}}}{{\binom{{2\,\ell + 2n}}{{2n}}\binom{{2m + 2n}}{{2n}}}}\,\frac{1}{\pi }\nonumber\\
&\qquad\qquad\qquad\times\left( {2O_{\ell + n}  + 2O_{m + n}  - H_{\ell + m + n - 1}  - 4\ln 2} \right)
\end{align}
and
\begin{align}\label{master3}
&\sum_{k = \ell}^\infty  {\frac{{\left( {\frac{1}{2}} \right)_k^2 }}{{k!^2 }}\,\frac{O_{k-\ell}}{{\prod_{j = 1}^{\,\ell} {\left( {2k - 2j + 1} \right)} \prod_{j = 1}^m {\left( {2k - 2j + 1} \right)} \prod_{j = 1}^n {\left(k + j\right)} }}}\nonumber\\
&\qquad  =-\sum_{k = 0}^{\ell-1}  {\frac{{\left( {\frac{1}{2}} \right)_k^2 }}{{k!^2 }}\,\frac{O_{\ell-k}}{{\prod_{j = 1}^{\,\ell} {\left( {2k - 2j + 1} \right)} \prod_{j = 1}^m {\left( {2k - 2j + 1} \right)} \prod_{j = 1}^n {\left(k + j\right)} }}}\nonumber\\
&\qquad\qquad\qquad +\frac{{\left( {\,\ell + m + n - 1} \right)!}}{{\left( {1/2} \right)_{\,\ell}^2 \left( {1/2} \right)_m^2 \left( {1/2} \right)_n^2 }}\,\frac{{( - 1)^{\,\ell + m} }}{{2^{\,\ell + m + 1} }}\,\frac{{\binom{{\,\ell + n}}{n}\binom{{m + n}}{n}}}{{\binom{{2\,\ell + 2n}}{{2n}}\binom{{2m + 2n}}{{2n}}}}\,\frac{1}{\pi }\nonumber\\
&\qquad\qquad\qquad\qquad\times\left( {2O_{\ell}  + 2O_{\ell + n}  - H_{\ell + m + n - 1}  - 2\ln 2} \right).
\end{align}

Identities~\eqref{master2} and~\eqref{master3} hold for non-negative integers $\ell$, $m$, and $n$ such that $\ell+m+n$ is a positive integer. At $\ell=1$, $m=1$, and $n=0$, Identity~\eqref{master2} gives the harmonic number complement of Forsyth's series~\eqref{forsyth} as
\begin{equation}
\sum_{k = 0}^\infty  {\frac{{\left( {\frac{1}{2}} \right)_k^2 }}{{k!^2 }}\,\frac{{H_k }}{{\left( {2k - 1} \right)^2 }}}  = \frac{4}{\pi }\,\left( {3 - 4\ln 2} \right),
\end{equation}
while $\ell=0$, $m=0$, and $n=1$ in~\eqref{master3} yields the odd harmonic number complement of Glaisher's series~\eqref{glaisher} as
\begin{equation}
\sum_{k = 0}^\infty  {\frac{{\left( {\frac{1}{2}} \right)_k^2 }}{{k!^2 }}\,\frac{{O_k }}{{k + 1}}}  = \frac{4}{\pi }\,\left( {1 - \ln 2} \right).
\end{equation}
Finally, we will derive identities involving products of harmonic numbers, namely,
\begin{align}\label{master4}
&\sum_{k = \ell}^\infty  {\frac{{\left( {\frac{1}{2}} \right)_k^2 }}{{k!^2 }}\,\frac{O_{k-\ell}H_{k+n}}{{\prod_{j = 1}^{\,\ell} {\left( {2k - 2j + 1} \right)} \prod_{j = 1}^m {\left( {2k - 2j + 1} \right)} \prod_{j = 1}^n {\left(k + j\right)} }}}\nonumber\\
&\qquad  =-\sum_{k = 0}^{\ell-1}  {\frac{{\left( {\frac{1}{2}} \right)_k^2 }}{{k!^2 }}\,\frac{O_{\ell-k}H_{k+n}}{{\prod_{j = 1}^{\,\ell} {\left( {2k - 2j + 1} \right)} \prod_{j = 1}^m {\left( {2k - 2j + 1} \right)} \prod_{j = 1}^n {\left(k + j\right)} }}}\nonumber\\ 
&\qquad\qquad+\frac{{\left( {\,\ell + m + n - 1} \right)!}}{{\left( {1/2} \right)_{\,\ell}^2 \left( {1/2} \right)_m^2 \left( {1/2} \right)_n^2 }}\,\frac{{( - 1)^{\,\ell + m} }}{{2^{\,\ell + m + 1} }}\,\frac{{\binom{{\,\ell + n}}{n}\binom{{m + n}}{n}}}{{\binom{{2\,\ell + 2n}}{{2n}}\binom{{2m + 2n}}{{2n}}}}\,\frac{1}{\pi }\nonumber\\
&\qquad\qquad\qquad\times\Big( {\left( {2O_{\ell + n}  + 2O_{m + n}  - H_{\ell + m + n - 1}  - 4\ln 2} \right)}\left( {2O_{\ell}  + 2O_{\ell + n}  - H_{\ell + m + n - 1}  - 2\ln 2} \right)\nonumber\\
&\qquad\qquad\qquad\qquad\qquad - 2\zeta (2) + 4O_{\ell + n}^{(2)}  - H_{\ell + m + n - 1}^{(2)} \Big),
\end{align}
with the particular case
\begin{equation}
\sum_{k = 0}^\infty  {\frac{{\left( {\frac{1}{2}} \right)_k^2 }}{{k!^2 }}\,\frac{{H_k O_k }}{{2k - 1}}}  =  - \frac{4}{\pi }\,\left( {2\ln 2 - 1} \right)\,\ln 2 + \frac{\pi }{3},
\end{equation}
and
\begin{align}\label{master5}
&\sum_{k = \max({\ell,m})}^\infty  {\frac{{\left( {\frac{1}{2}} \right)_k^2 }}{{k!^2 }}\,\frac{O_{k-\ell}O_{k-m}}{{\prod_{j = 1}^{\,\ell} {\left( {2k - 2j + 1} \right)} \prod_{j = 1}^m {\left( {2k - 2j + 1} \right)} \prod_{j = 1}^n {\left(k + j\right)} }}}\nonumber\\
&\qquad  =-\sum_{k = 0}^{\max({\ell,m})-1}  {\frac{{\left( {\frac{1}{2}} \right)_k^2 }}{{k!^2 }}\,\frac{O_{\ell-k}O_{m-k}}{{\prod_{j = 1}^{\,\ell} {\left( {2k - 2j + 1} \right)} \prod_{j = 1}^m {\left( {2k - 2j + 1} \right)} \prod_{j = 1}^n {\left(k + j\right)} }}}\nonumber\\ 
&\qquad\qquad+\frac{{\left( {\,\ell + m + n - 1} \right)!}}{{\left( {1/2} \right)_{\,\ell}^2 \left( {1/2} \right)_m^2 \left( {1/2} \right)_n^2 }}\,\frac{{( - 1)^{\,\ell + m} }}{{2^{\,\ell + m + 2} }}\,\frac{{\binom{{\,\ell + n}}{n}\binom{{m + n}}{n}}}{{\binom{{2\,\ell + 2n}}{{2n}}\binom{{2m + 2n}}{{2n}}}}\,\frac{1}{\pi }\nonumber\\
&\qquad\qquad\qquad\times\Big( {\left( {2O_m  + 2O_{m + n}  - H_{\ell + m + n - 1}  - 2\ln 2} \right)}\left( {2O_{\ell}  + 2O_{\ell + n}  - H_{\ell + m + n - 1}  - 2\ln 2} \right)\nonumber\\
&\qquad\qquad\qquad\qquad\qquad + \zeta (2) - H_{\ell + m + n - 1}^{(2)} \Big),
\end{align}
with the special value
\begin{equation}
\sum_{k = 0}^\infty  {\frac{{\left( {\frac{1}{2}} \right)_k^2 }}{{k!^2 }}\,\frac{{O_k^2 }}{{k + 1}}}  = \frac{4}{\pi }\left( {1 - \ln 2} \right)^2  + \frac{\pi }{6}.
\end{equation}

In~\eqref{master4} and~\eqref{master5}, $H_j^{(2)}$ and $O_j^{(2)}$ are the second order $j^{th}$ harmonic and odd harmonic number, respectively. Here, and throughout this paper, $\zeta(s)$ denotes the Riemann zeta function. 

Identities~\eqref{master4} and~\eqref{master5} are valid for non-negative integers $\ell$, $m$, and $n$ such that $\ell+m+n$ is a positive integer.

\section{Preliminary results}
\begin{lemma}[{\cite[Lemma 2]{adegoke0}}]\label{zcr19hm}
If $n$ is a non-negative integer, then
\begin{equation}\label{mq7ny8m}
\left( { - n + \frac{1}{2}} \right)_k  = ( - 1)^n \frac{{2^n \left( {\frac{1}{2}} \right)_n \left( {\frac{1}{2}} \right)_k }}{{\prod_{j = 1}^n {\left( {2k - 2j + 1} \right)} }}.
\end{equation}
\end{lemma}
\begin{lemma}\label{feoqz9e}
If $\ell$, $m$, and $n$ are non-negative integers whose sum is a positive integer, then
\begin{equation}\label{mxp1ake}
\left( {\frac{1}{2} + \ell + n} \right)_{m - \frac{1}{2}} \left( {\frac{1}{2} + m + n} \right)_{\frac{1}{2} - m}  = \frac{{n!\left( {\ell + m + n - 1} \right)!}}{{\left( {\frac{1}{2}} \right)_n^2 \left( {\frac{1}{2} + n} \right)_{\ell} \left( {\frac{1}{2} + n} \right)_m }}\,\frac{1}{\pi }.
\end{equation}
Observe the symmetry property
\begin{equation}\label{o71grqs}
\left( {\frac{1}{2} + \ell + n} \right)_{m - \frac{1}{2}} \left( {\frac{1}{2} + m + n} \right)_{\frac{1}{2} - m}=\left( {\frac{1}{2} + m + n} \right)_{\ell - \frac{1}{2}} \left( {\frac{1}{2} + \ell + n} \right)_{\frac{1}{2} - \ell}.
\end{equation}
\end{lemma}

Let $\Gamma(z)$ be Euler's Gamma function. Generalized binomial coefficients are defined for complex numbers $r$ and $s$ by
\begin{equation}\label{y89d722}
\binom rs= \frac{{\Gamma (r + 1)}}{{\Gamma (s + 1)\Gamma (r - s + 1)}}.
\end{equation}
The {\it extended} Pochhammer relation $(x)_y$ (see~\cite[Definition 1.1.]{liu12}) is defined for complex numbers $x$ and $y$ that are not negative integers by
\begin{equation}\label{ex_pochhammer}
\left( x \right)_y  = \frac{{\Gamma \left( {x + y} \right)}}{{\Gamma \left( x \right)}} = \binom{{x + y - 1}}{y}\,\Gamma \left( {y + 1} \right).
\end{equation}
If, in particular,  $y=n$ is a non-negative integer, then we have the Pochhammer symbol
\begin{equation}
\left( x \right)_n  = \prod_{k = 0}^{n - 1} {\left(x + k\right)}  = \binom{{x + n - 1}}{n}\,n!
\end{equation}

For non-negative integers $n$ and positive integers (order) $m$, the standard generalized harmonic numbers and generalized odd harmonic numbers are defined by 
\begin{equation*}
H_n^{(m)} = \sum_{k=1}^n \frac{1}{k^m} \text{ and } O_n^{(m)} = \sum_{k=1}^n \frac{1}{(2k-1)^m},
\end{equation*}
respectively, with $H_0^{(m)} = O_0^{(m)} = 0$. We further adopt the notation: \mbox{$H_n^{(1)} =H_n$}, \mbox{$O_n^{(1)} =O_n$}.
 
By extending these definitions to fractional arguments via the polygamma functions, a direct algebraic connection emerges at half-integer coordinates. Specifically, the standard harmonic number of order $m+1$ at a half-integer argument relates explicitly to the odd harmonic number of the same order via the identity
\begin{equation}\label{g9e4dpz}
H_{n-1/2}^{(m+1)} - H_{-1/2}^{(m+1)}= 2^{m+1} O_n^{(m+1)}.
\end{equation}

This evaluation allows half-integer evaluations of the polygamma functions to be converted seamlessly into discrete, integer-indexed sums.

When $m=0$, the relationship between harmonic numbers, odd harmonic numbers and the digamma function establishes that \mbox{$H_{-1/2}=-2\ln 2$}, so that in this case~\eqref{g9e4dpz} reduces to
\begin{equation}\label{nerd5x1}
H_{n-1/2}=2O_n-2\ln 2.
\end{equation}
\begin{lemma}\label{ygxw07j}
We have
\begin{equation}\label{ty98d9m}
\frac{\partial}{{\partial x}}\left( x \right)_y  = \left( x \right)_y\left(\psi(x+y)-\psi(x)\right)= \left( x \right)_y \left( {H_{x + y - 1}  - H_{x - 1} } \right),
\end{equation}
where $\psi(x)$ is the digamma function.
\end{lemma}

\begin{lemma}[{\cite{xu17}}]\label{qlurtgo}
If $\alpha$ and $u$ are complex numbers such that $\alpha+u$ lie within the domain of analyticity, then
\begin{equation*}
\frac{{\partial ^m }}{{\partial \alpha ^m }}H_{u + \alpha }^{(m)}  = m\left( {\zeta (m + 1) - H_{u + \alpha }^{(m + 1)} } \right).
\end{equation*}
In particular,
\begin{equation*}
\frac{\partial }{{\partial \alpha }}H_{u + \alpha }  = \zeta (2) - H_{u + \alpha }^{(2)} .
\end{equation*}
\end{lemma}

\section{The master identity and associated identities}
\begin{theorem}
If $\ell$, $m$, and $n$ are non-negative integers whose sum is a positive integer, then
\begin{align*}
&\sum_{k = 0}^\infty  {\frac{{\left( {\frac{1}{2}} \right)_k^2 }}{{k!^2 }}\,\frac{1}{{\prod_{j = 1}^{\,\ell} {\left( {2k - 2j + 1} \right)} \prod_{j = 1}^m {\left( {2k - 2j + 1} \right)} \prod_{j = 1}^n {\left(k + j\right)} }}}\nonumber\\
&\qquad  = \frac{{\left( {\,\ell + m + n - 1} \right)!}}{{\left( {1/2} \right)_{\,\ell}^2 \left( {1/2} \right)_m^2 \left( {1/2} \right)_n^2 }}\,\frac{{( - 1)^{\,\ell + m} }}{{2^{\,\ell + m} }}\,\frac{{\binom{{\,\ell + n}}{n}\binom{{m + n}}{n}}}{{\binom{{2\,\ell + 2n}}{{2n}}\binom{{2m + 2n}}{{2n}}}}\,\frac{1}{\pi }\tag{\ref{master1}}.
\end{align*}
\end{theorem}
\begin{proof}
Setting $a=1/2-p$, $b=1/2-q$, and $c=r+1$ in Gauss's hypergeometric summation formula
\begin{equation*}
\sum_{k = 0}^\infty  {\frac{{\left( a \right)_k \left( b \right)_k }}{{k!\left( c \right)_k }}}  = \left( {c - a} \right)_a \left( {c - b} \right)_{ - a} ,\quad\Re\left( {c - a - b} \right) > 0,
\end{equation*}
gives
\begin{equation}\label{h0piaxg}
\sum_{k = 0}^\infty  {\frac{{\left( {1/2 - p} \right)_k \left( {1/2 - q} \right)_k }}{{k!\left( {1 + r} \right)_k }}}  = \left( {1/2 + p + r} \right)_{q - 1/2} \left( {1/2 + q + r} \right)_{1/2 - q} ,
\end{equation}
for $\Re(p+q+r)>0$, from which~\eqref{master1} follows upon setting $p=\ell$, $q=m$, and $r=n$ and using Lemmata~\ref{zcr19hm} and~\ref{feoqz9e} and the fact that
\begin{equation}
\left( {1 + n} \right)_k  = \frac{{k!}}{{n!}}\prod_{j = 1}^n {\left( {k + j} \right)} .
\end{equation}
\end{proof}
In the remaining part of this section we will derive associated identities to the master identity, involving the standard and odd harmonic numbers and their products. These will be obtained via differentiation with respect to the parameters $p$, $q$, and $r$. Termwise differentiation is justified locally in the parameter domain by the convergence of Gauss's hypergeometric series and its parameter derivatives.
\begin{theorem}
If $\ell$, $m$, and $n$ are non-negative integers whose sum is a positive integer, then
\begin{align*}
&\sum_{k = 0}^\infty  {\frac{{\left( {\frac{1}{2}} \right)_k^2 }}{{k!^2 }}\,\frac{H_{k+n}}{{\prod_{j = 1}^{\,\ell} {\left( {2k - 2j + 1} \right)} \prod_{j = 1}^m {\left( {2k - 2j + 1} \right)} \prod_{j = 1}^n {\left(k + j\right)} }}}\nonumber\\
&\qquad  = \frac{{\left( {\,\ell + m + n - 1} \right)!}}{{\left( {1/2} \right)_{\,\ell}^2 \left( {1/2} \right)_m^2 \left( {1/2} \right)_n^2 }}\,\frac{{( - 1)^{\,\ell + m} }}{{2^{\,\ell + m} }}\,\frac{{\binom{{\,\ell + n}}{n}\binom{{m + n}}{n}}}{{\binom{{2\,\ell + 2n}}{{2n}}\binom{{2m + 2n}}{{2n}}}}\,\frac{1}{\pi }\nonumber\\
&\qquad\qquad\qquad\times\left( {2O_{\ell + n}  + 2O_{m + n}  - H_{\ell + m + n - 1}  - 4\ln 2} \right)\tag{\ref{master2}},
\end{align*}
and
\begin{align*}
&\sum_{k = \ell}^\infty  {\frac{{\left( {\frac{1}{2}} \right)_k^2 }}{{k!^2 }}\,\frac{O_{k-\ell}}{{\prod_{j = 1}^{\,\ell} {\left( {2k - 2j + 1} \right)} \prod_{j = 1}^m {\left( {2k - 2j + 1} \right)} \prod_{j = 1}^n {\left(k + j\right)} }}}\nonumber\\
&\qquad  =-\sum_{k = 0}^{\ell-1}  {\frac{{\left( {\frac{1}{2}} \right)_k^2 }}{{k!^2 }}\,\frac{O_{\ell-k}}{{\prod_{j = 1}^{\,\ell} {\left( {2k - 2j + 1} \right)} \prod_{j = 1}^m {\left( {2k - 2j + 1} \right)} \prod_{j = 1}^n {\left(k + j\right)} }}}\nonumber\\
&\qquad\qquad\qquad +\frac{{\left( {\,\ell + m + n - 1} \right)!}}{{\left( {1/2} \right)_{\,\ell}^2 \left( {1/2} \right)_m^2 \left( {1/2} \right)_n^2 }}\,\frac{{( - 1)^{\,\ell + m} }}{{2^{\,\ell + m + 1} }}\,\frac{{\binom{{\,\ell + n}}{n}\binom{{m + n}}{n}}}{{\binom{{2\,\ell + 2n}}{{2n}}\binom{{2m + 2n}}{{2n}}}}\,\frac{1}{\pi }\nonumber\\
&\qquad\qquad\qquad\qquad\times\left( {2O_{\ell}  + 2O_{\ell + n}  - H_{\ell + m + n - 1}  - 2\ln 2} \right)\tag{\ref{master3}}.
\end{align*}
\end{theorem}
\begin{proof}
Differentiating~\eqref{h0piaxg} with respect to $r$ using Lemma~\ref{ygxw07j} gives
\begin{align}\label{aa09rq5}
\sum_{k = 0}^\infty  {\frac{{\left( {1/2 - p} \right)_k \left( {1/2 - q} \right)_k }}{{k!\left( {1 + r} \right)_k }}\,H_{k + r} }  &= \left( {1/2 + p + r} \right)_{q - 1/2} \left( {1/2 + q + r} \right)_{1/2 - q}\nonumber\\
&\qquad\qquad\times \left( {H_{r + p - 1/2}  + H_{r + q - 1/2}  - H_{r + p + q - 1} } \right),
\end{align}
after simplifying the algebra. Setting $p=\ell$, $q=m$, and $r=n$ and applying Lemmata~\ref{zcr19hm} and~\ref{feoqz9e} and Equation~\eqref{nerd5x1} yields~\eqref{master2}. To prove~\eqref{master3}, differentiate~\eqref{h0piaxg} with respect to~$p$ to obtain
\begin{align*}
\sum_{k = 0}^\infty  {\frac{{\left( {1/2 - p} \right)_k \left( {1/2 - q} \right)_k }}{{k!\left( {1 + r} \right)_k }}\left( {H_{k - p - 1/2}  - H_{ - p - 1/2} } \right)}  &= \left( {1/2 + p + r} \right)_{q - 1/2} \left( {1/2 + q + r} \right)_{1/2 - q}\nonumber\\
&\qquad\qquad \times\left( {H_{r + p - 1/2}  - H_{r + p + q - 1} } \right),
\end{align*} 
which, by Equation~\eqref{nerd5x1}, simplifies further to
\begin{align}\label{uy817hg}
\sum_{k = 0}^\infty  {\frac{{\left( {1/2 - p} \right)_k \left( {1/2 - q} \right)_k }}{{k!\left( {1 + r} \right)_k }}\,2O_{k-p}}  &= \left( {1/2 + p + r} \right)_{q - 1/2} \left( {1/2 + q + r} \right)_{1/2 - q}\nonumber\\
&\qquad\qquad \times\left( {2O_p + 2O_{r + p}  - H_{r + p + q - 1}-2\ln 2 } \right).
\end{align}
Identity~\eqref{master3} now follows by plugging $p=\ell$, $q=m$, and $r=n$ and applying Lemmata~\ref{zcr19hm} and~\ref{feoqz9e}. Note that the final form of~\eqref{master3} follows from the fact that $O_{-j}=O_j$.
\end{proof}
\begin{lemma}\label{kvwevuq}
If $p$, $q$, and $r$ are complex numbers such that the arguments lie within the domain of analyticity, then
\begin{equation}\label{myb5sq3}
\frac{\partial }{{\partial p}}\left( {H_{r + p - 1/2}  - H_{r + p + q - 1} } \right) = H_{r + p + q - 1}^{(2)}  + 2\zeta (2) - 4O_{r + p}^{(2)} .
\end{equation}
\end{lemma}
\begin{proof}
Using Lemma~\ref{qlurtgo}, we have
\begin{equation}\label{b7fqx56}
\frac{\partial }{{\partial p}}\left( {H_{r + p - 1/2}  - H_{r + p + q - 1} } \right) = H_{r + p + q - 1}^{(2)}  - H_{r + p - 1/2}^{(2)} .
\end{equation}
Using $m=1$ and $n=r+p$ in~\eqref{g9e4dpz} gives
\begin{equation}\label{v3d2plw}
H_{r + p - 1/2}^{(2)}  = H_{ - 1/2}^{(2)}  + 4O_{r + p}^{(2)} .
\end{equation}
But, since (see~\cite{sofo16})
\begin{equation*}
H_{1/2}^{(2)}  = 4  - 2\zeta (2) ,
\end{equation*}
we have
\begin{equation*}
H_{ - 1/2}^{(2)}  =  - 2\zeta (2),
\end{equation*}
which, substituted in~\eqref{v3d2plw} yields
\begin{equation}\label{q1uzle0}
H_{r + p - 1/2}^{(2)}  =  - 2\zeta (2) + 4O_{r + p}^{(2)} .
\end{equation}
Substituting~\eqref{q1uzle0} in~\eqref{b7fqx56} gives~\eqref{myb5sq3}.
\end{proof}
\begin{theorem}
If $\ell$, $m$, and $n$ are non-negative integers whose sum is a positive integer, then
\begin{align*}
&\sum_{k = \ell}^\infty  {\frac{{\left( {\frac{1}{2}} \right)_k^2 }}{{k!^2 }}\,\frac{O_{k-\ell}H_{k+n}}{{\prod_{j = 1}^{\,\ell} {\left( {2k - 2j + 1} \right)} \prod_{j = 1}^m {\left( {2k - 2j + 1} \right)} \prod_{j = 1}^n {\left(k + j\right)} }}}\nonumber\\
&\qquad  =-\sum_{k = 0}^{\ell-1}  {\frac{{\left( {\frac{1}{2}} \right)_k^2 }}{{k!^2 }}\,\frac{O_{\ell-k}H_{k+n}}{{\prod_{j = 1}^{\,\ell} {\left( {2k - 2j + 1} \right)} \prod_{j = 1}^m {\left( {2k - 2j + 1} \right)} \prod_{j = 1}^n {\left(k + j\right)} }}}\nonumber\\ 
&\qquad\qquad+\frac{{\left( {\,\ell + m + n - 1} \right)!}}{{\left( {1/2} \right)_{\,\ell}^2 \left( {1/2} \right)_m^2 \left( {1/2} \right)_n^2 }}\,\frac{{( - 1)^{\,\ell + m} }}{{2^{\,\ell + m + 1} }}\,\frac{{\binom{{\,\ell + n}}{n}\binom{{m + n}}{n}}}{{\binom{{2\,\ell + 2n}}{{2n}}\binom{{2m + 2n}}{{2n}}}}\,\frac{1}{\pi }\nonumber\\
&\qquad\qquad\qquad\times\Big( {\left( {2O_{\ell + n}  + 2O_{m + n}  - H_{\ell + m + n - 1}  - 4\ln 2} \right)}\left( {2O_{\ell}  + 2O_{\ell + n}  - H_{\ell + m + n - 1}  - 2\ln 2} \right)\nonumber\\
&\qquad\qquad\qquad\qquad\qquad - 2\zeta (2) + 4O_{\ell + n}^{(2)}  - H_{\ell + m + n - 1}^{(2)} \Big)\tag{\ref{master4}}
\end{align*}
and
\begin{align*}
&\sum_{k = \max({\ell,m})}^\infty  {\frac{{\left( {\frac{1}{2}} \right)_k^2 }}{{k!^2 }}\,\frac{O_{k-\ell}O_{k-m}}{{\prod_{j = 1}^{\,\ell} {\left( {2k - 2j + 1} \right)} \prod_{j = 1}^m {\left( {2k - 2j + 1} \right)} \prod_{j = 1}^n {\left(k + j\right)} }}}\nonumber\\
&\qquad  =-\sum_{k = 0}^{\max({\ell,m})-1}  {\frac{{\left( {\frac{1}{2}} \right)_k^2 }}{{k!^2 }}\,\frac{O_{\ell-k}O_{m-k}}{{\prod_{j = 1}^{\,\ell} {\left( {2k - 2j + 1} \right)} \prod_{j = 1}^m {\left( {2k - 2j + 1} \right)} \prod_{j = 1}^n {\left(k + j\right)} }}}\nonumber\\ 
&\qquad\qquad+\frac{{\left( {\,\ell + m + n - 1} \right)!}}{{\left( {1/2} \right)_{\,\ell}^2 \left( {1/2} \right)_m^2 \left( {1/2} \right)_n^2 }}\,\frac{{( - 1)^{\,\ell + m} }}{{2^{\,\ell + m + 2} }}\,\frac{{\binom{{\,\ell + n}}{n}\binom{{m + n}}{n}}}{{\binom{{2\,\ell + 2n}}{{2n}}\binom{{2m + 2n}}{{2n}}}}\,\frac{1}{\pi }\nonumber\\
&\qquad\qquad\qquad\times\Big( {\left( {2O_m  + 2O_{m + n}  - H_{\ell + m + n - 1}  - 2\ln 2} \right)}\left( {2O_{\ell}  + 2O_{\ell + n}  - H_{\ell + m + n - 1}  - 2\ln 2} \right)\nonumber\\
&\qquad\qquad\qquad\qquad\qquad + \zeta (2) - H_{\ell + m + n - 1}^{(2)} \Big)\tag{\ref{master5}}.
\end{align*}
\end{theorem}
\begin{proof}
Differentiate~\eqref{aa09rq5} with respect to $p$ and simplify, with the help of Equation~\eqref{nerd5x1} and Lemma~\ref{kvwevuq}, to get
\begin{align*}
&\sum_{k = 0}^\infty  {\frac{{\left( {1/2 - p} \right)_k \left( {1/2 - q} \right)_k }}{{k!\left( {1 + r} \right)_k }}\,2O_{k - p} H_{k + r} }\\
&\qquad  = \left( {1/2 + p + r} \right)_{q - 1/2} \left( {1/2 + q + r} \right)_{1/2 - q}\Big(\\ 
&\qquad\qquad {\left( {2O_{r + p}  + 2O_{r + q}  - H_{r + p + q - 1}  - 4\ln 2} \right)\left( {2O_p  + 2O_{r + p}  - H_{r + p + q - 1}  - 2\ln 2} \right)}\\
&\qquad\qquad\qquad { - 2\zeta (2) + 4O_{r + p}^{(2)}  - H_{r + p + q - 1}^{(2)} } \Big).
\end{align*}
Now substitute $p=\ell$, $q=m$, and $r=n$ and apply Lemmata~\ref{zcr19hm} and~\ref{feoqz9e}. This completes the proof of~\eqref{master4}. To prove~\eqref{master5}, employ the symmetry property~\eqref{o71grqs} to write~\eqref{uy817hg} as
\begin{align*}
\sum_{k = 0}^\infty  {\frac{{\left( {1/2 - p} \right)_k \left( {1/2 - q} \right)_k }}{{k!\left( {1 + r} \right)_k }}\,2O_{k-p}}  &= \left( {1/2 + q + r} \right)_{p - 1/2} \left( {1/2 + p + r} \right)_{1/2 - p}\nonumber\\
&\qquad\qquad \times\left( {2O_p + 2O_{r + p}  - H_{r + p + q - 1}-2\ln 2 } \right).
\end{align*}
Differentiate with respect to $q$, using Lemmata~\ref{ygxw07j} and~\ref{qlurtgo}, to obtain, after some algebra
\begin{align*}
&\sum_{k = 0}^\infty  {\frac{{\left( {1/2 - p} \right)_k \left( {1/2 - q} \right)_k }}{{k!\left( {1 + r} \right)_k }}\,4O_{k - p} O_{k - q} }\\
&\qquad  = \left( {1/2 + p + r} \right)_{q - 1/2} \left( {1/2 + q + r} \right)_{1/2 - q}\Big(\\ 
&\qquad\qquad {\left( {2O_{r + p}  + 2O_p  - H_{r + p + q - 1}  - 2\ln 2} \right)\left( {2O_q  + 2O_{r + q}  - H_{r + p + q - 1}  - 2\ln 2} \right)}\\
&\qquad\qquad\qquad { +\zeta (2)  - H_{r + p + q - 1}^{(2)} } \Big).
\end{align*}
Set $p=\ell$, $q=m$, and $r=n$ and apply Lemmata~\ref{zcr19hm} and~\ref{feoqz9e} to obtain~\eqref{master5}.
\end{proof}

\section{Concluding Remarks}
 In this work, we have established a robust three-parameter master identity that unifies and contextualizes the landscape of Forsyth-Glaisher-type series for $1/\pi$. By anchoring our proof in Gauss's classical hypergeometric summation formula, we have shown that historically isolated formulas including Glaisher's original 1905 series are merely specific structural manifestations of a broader algebraic truth.

The primary utility of this master identity lies in its generative and interpolative power. It not only simplifies the recovery of known literature but opens an avenue to produce an infinite family of entirely new $1/\pi$ expansions. Furthermore, our exploration into the parameter space via differentiation has yielded elegant, unexpected connection points to both standard and odd harmonic numbers.

Several promising directions for future research emerge from this framework:
\begin{itemize}
    \item \textbf{Higher-order Derivatives:} Investigating the consequence of taking higher-order derivatives with respect to the parameters may unlock more intricate series involving generalized harmonic numbers of higher depths.
    \item \textbf{Generalization to Higher Hypergeometric Functions:} A natural next step is to explore whether a similar multi-parameter master framework can be extended to $_3F_2$ or $_4F_3$ series, potentially unifying Ramanujan's higher-order $1/\pi$ series or even $1/\pi^2$ expansions.
    \item \textbf{p-adic Analogues:} Exploring the p-adic and q-analogues of this master identity could unveil deep arithmetic properties and supercongruences akin to those found in modern number theory.
\end{itemize}
Ultimately, this unified view underscores that classical hypergeometric tools still possess immense untapped potential for organizing and expanding our understanding of mathematical constants.

\end{document}